\documentstyle[12pt]{article}

\font\iBbb = msbm9
\font\Bbb = msbm10 scaled\magstep1

\font\ttBbb = msbm10 scaled\magstep3
\font\tttBbb = msbm10 scaled\magstep4
\font\goth = eufm10 scaled\magstep1

\font\ttgoth = eufm10 scaled\magstep3
\font\tttgoth = eufm10 scaled\magstep4

\def\Stovicek{{\v{S}\v{t}ov\'\i\v{c}ek\ }}

\def\iR{\mbox{\iBbb R}}
\def\R{\mbox{\Bbb R}}

\def\ttR{\mbox{\ttBbb R}}
\def\tttR{\mbox{\tttBbb R}}
\def\C{\mbox{\Bbb C}}

\def\Z{\mbox{\Bbb Z}}
\def\N{\mbox{\Bbb N}}

\def\sl{\mbox{\goth sl}} 
 
\def\ttsl{\mbox{\ttgoth sl}} 
\def\tttsl{\mbox{\tttgoth sl}} 
\def\g{\mbox{\goth g}}
\def\gg{{\goth g}}
\def\su{\mbox{\goth su}} 
\def\so{\mbox{\goth so}} 

\def\UU{{\cal U}}
\def\CC{{\cal C}}

\def\AA{{\cal A}}

\def\Tr{\mbox{Tr}}
\def\sgn{\mbox{sgn}}
\def\Im{\mbox{Im}}

\def\Usl{{\cal U}_q(\sl(2,\R))}
\def\slR{\sl(2,\R)}
\def\UslC{{\cal U}_q(\sl(2,\C))}
\def\ASLC{{\cal A}_q(SL(2,\C))}

\normalbaselines 
\begin{document} 
\begin{center}
\vspace*{1.0cm}

{\LARGE{\bf Discrete series of representations for 
\protect\( \UU _{q}(\tttsl (2,\tttR ))\protect \)}} 

\vskip 1.5cm

{\large {\bf Pavel \Stovicek }} 

\vskip 0.5 cm 

Department of Mathematics, Faculty of Nuclear Science \\ 
Czech Technical University \\ 
Trojanova 13, 120 00 Prague\\ 
Czech Republic

\end{center}

\vspace{1 cm}

\begin{abstract}
A possible generalization of the method of orbits to \( SL_{q}(2,\R ) \) 
is discussed. 
\end{abstract}

\vspace{1 cm} 

\section{Introduction}

There are not so many results concerning non-compact real quantum groups. 
Naturally a major attention has been paid to groups with some physical
interpretation (see, for example, \cite{DM}). Even in the case of real groups
of low dimension the treatment seems to be quite complicated \cite{Woronowicz}.
Recently the quantum group \( SL_{q}(2,\R ) \) attracted some attention; particular
efforts have been made to describe and classify its quantum homogeneous spaces
\cite{Tarlini}. In this contribution we focus on a possible generalization 
of the method of orbits to \(\UU _{q}(\sl(2,\R ))\). 

In the paper \cite{PS} there was described a construction of representations
for a Hopf algebra \( \UU  \) which may be considered
as a generalization of the method of orbits due to Kostant and Kirillov. In
fact, the construction concentrates on just one step of the method of orbits.
It is already supposed that there is given a left \( \UU  \)-module 
algebra \( \CC  \)
with the action denoted by \( \xi  \) and fulfilling two conditions: 
$\xi (X)\cdot 1=\varepsilon (X)\, 1$ and Leibniz rule 
\begin{equation}
\xi (X)\cdot (fg)=(\xi (X_{(1)})\cdot f)(\xi (X_{(2)})\cdot g),\quad 
\forall X\in \UU ,\, \forall f,g\in \CC .\label{Leibniz} 
\end{equation}
We use Sweedler's notation: \( \Delta X=X_{(1)}\otimes X_{(2)} \) where 
\( \Delta  \) denotes the comultiplication. 

Classically, \( \UU  \) is the universal enveloping algebra of a Lie algebra
\gg, \( \CC  \) is an algebra of functions living on a coadjoint orbit, or
a reduction of such an algebra via a polarization, and \( \xi  \) is nothing
but the infinitesimal action. So \( \xi (X) \), \( X\in \g  \), are vector
fields. A critical step of the method of orbits is a replacement of vector fields
\( \xi (X) \) by first order differential operators of the form 
\( \nabla (\xi (X))+2\pi \imath \lambda (X) \)
where \( \nabla  \) is an appropriate covariant derivative and \( \lambda (X) \)
is a function depending linearly on \( X \). And this is where the construction
attempts a generalization applicable also to quantum groups. It results in  
a modification of the left action \( \xi  \). The new action will
be denoted simply by a dot. 

The construction is also known in somewhat different form as the twisted 
adjoint action \cite{Joseph}.

%

The simplest example is \( \UslC  \) (\( q \) is not a root of unity) with
the generators \( K,\, K^{-1},\, E \) and \( F \), and the defining relations
\begin{eqnarray}
 & K\, K^{-1}=K^{-1}K=1,\, \, KE=q\, EK,\, \, KF=q^{-1}FK, & \nonumber \\
 & [E,F]=\frac{1}{q-q^{-1}}(K^{2}-K^{-2}). & 
\end{eqnarray}
The comultiplication, the counite and the antipode are the usual ones.  

Let \( \AA =\ASLC  \) be the dual 
Hopf algebra of quantum functions on \( SL(2,\C ) \)
with the generators \( a,\, b,\, c,\, d \), and set \( \CC =\C [z] \). Then 
\begin{equation}
L:\CC \to \AA \otimes \CC, \, 
L(z)=(c\otimes 1+d\otimes z)(a\otimes 1+b\otimes z)^{-1}, 
\end{equation}
is a left coaction (defined in fact rather formally). 
Using the standard pairing between
\( \UslC  \) and \( \ASLC  \) one introduces a right action, and  
with the aid of the antipode
\( S \) one can pass to a left action, 
\begin{equation}
\xi (x)\cdot f=\langle SX,f_{(1)}\rangle \, f_{(2)},\quad 
\forall X\in \UU ,\, \forall f\in \CC .
\end{equation}

After some rescaling of the complex variable \( z \) one arrives at the action
\begin{equation}
\xi (K)\cdot z^{j}=q^{j}z^{j},\, \, 
\xi (E)\cdot z^{j}=[j]_{q}\,z^{j+1},\, \, 
\xi (F)\cdot z^{j}=-[j]_{q}\,z^{j-1},
\end{equation}
where 
\begin{equation}
[x]_{q}=\frac{q^{x}-q^{-x}}{q-q^{-1}}.
\end{equation}
The mentioned construction is applicable to this case and yields a one-parameter
family of modified actions 
\begin{equation}
\label{maction}
K\cdot z^{j}=q^{(2j-\sigma )/2}z^{j},\, \, 
E\cdot z^{j}=q^{-\sigma /2}[j-\sigma ]_{q}\,z^{j+1},\, \, 
F\cdot z^{j}=-q^{\sigma /2}[j]_{q}\,z^{j-1}.
\end{equation}
The action can be restricted to the two real forms of \( \UslC  \), namely 
\( \UU _{q}(\su (2)) \)
and \( \Usl  \). The latter one is determined
by the involution 
\begin{equation}
\label{slinvolution}
K^{\ast }=K,\, \, E^{\ast }=-q^{-1}\, E,\, \, F^{\ast }=-q\, F,
\end{equation}
with \( q \) being a complex unite (\( q^{\ast }=q^{-1} \)).

As shown in \cite{PS}, the construction is applicable to
any compact simple Lie group from the four principal series and yields all irreducible
representations. In the present paper we concentrate, however, on the non-compact
real form \( \Usl  \).

\section{Classical method of orbits for \protect\( \ttsl (2,\ttR )\protect \)}

Let us identify \( \slR  \) with \( \R ^{3} \): 
\begin{equation}
\R ^{3}\ni (x_{1},x_{2},x_{3})\mapsto \left( \begin{array}{cc}
x_{3} & x_{1}+x_{2}\\
-x_{1}+x_{2} & -x_{3}
\end{array}\right) \in \slR .
\end{equation}
Relating to every matrix \( X\in \slR  \) a functional \( f_{X} \),
$\langle f_{X},Y\rangle =\Tr \, X^{t}Y$,
we identify \( \slR  \) with its dual. The coadjoint orbits are either the
origin or the quadrics 
\begin{equation}
x_{1}^{\, 2}-x_{2}^{\, 2}-x_{3}^{\, 2}=c=\textrm{const}.
\end{equation}
The orbits contributing to the discrete series are characterized by 
$c=k^{2},\, \, k>0$. 
We choose the leaf with \( \sgn \, x_{1}>0 \). Then it is natural to use 
\( (x_{2},x_{3}) \)
as coordinates on the orbit, with \( x_{1}=(k^{2}+x_{2}^{\, 2}+x_{3}^{\, 2})^{1/2} \). 

However it is more convenient to introduce complex coordinates \( (z,\bar{z}) \)
as follows: 
\begin{equation}
z=\frac{x_{1}+x_{2}}{x_{3}-\imath k}.
\end{equation}
Then the orbit is identified with the upper complex half-plane, \( \Im \, z>0 \),
and the coadjoint action takes the form 
\begin{equation}
g\cdot z=\frac{c+dz}{a+bz}\quad \textrm{where }g=\left( \begin{array}{cc}
a & b\\
c & d
\end{array}\right) \in SL(2,\R ).
\end{equation}
The infinitesimal action reads ($E,\ F,\ H$ is the usual basis) 
\begin{equation}
\label{infaction}
\xi _{E}=z^{2}\, \partial _{z}+\textrm{c}.\textrm{c}.,\, \, 
\xi _{F}=-\partial _{z}+\textrm{c}.\textrm{c}.,\, \, 
\xi _{H}=2z\, \partial _{z}+\textrm{c}.\textrm{c}.
\end{equation}

There exist two polarizations which are nothing but the mutually conjugate complex
structures \( \partial _{\overline{z}} \) and \( \partial _{z} \). We choose
the former one and so the constructed representation should act in a space
of holomorphic functions on the upper complex half-plane. 

The method of orbits leads to the following modification of the infinitesimal
action (\ref{infaction}):
\begin{equation}
\varrho (E)=z^{2}\, \partial _{z}+4\pi k,\, \, 
\varrho (F)=-\partial _{z},\, \, \varrho (H)=2z\, \partial _{z}+4\pi k.
\end{equation}
This representation of the Lie algebra can be integrated to a representation
of the Lie group, and the result reads 
\begin{equation}
\varrho (g)\cdot f(z)=(d-bz)^{-4\pi k}f\left( \frac{-c+az}{d-bz}\right) \quad 
\textrm{where }g=\left( \begin{array}{cc}
a & b\\
c & d
\end{array}\right) .
\end{equation}
The representation \( \varrho  \) depends on \( g \) continuously and unambiguously
if and only if 
\begin{equation}
4\pi k=n\in \N .
\end{equation}
This constraint is caused by the nontrivial topology of \( SL(2,\R ) \). Moreover,
the representation \( \varrho  \) is unitary with respect to the \( L^{2} \)-norm

\begin{equation}
\Vert f\Vert ^{2}=\int _{\iR \times \iR _{+}}|f(u,v)|^{2}v^{n-2}\, du\, dv
\end{equation}
where \( z=u+\imath v \). Because of the divergence at \( v=0 \) we have to
impose another constraint, namely $n\geq 2$. 

The restriction of \( \varrho  \) to the compact subgroup 
\( SO(2)\subset SL(2,\R ) \)
can be diagonalized and the eigen-functions form an orthogonal basis in the
carrier Hilbert space. The element \( E-F \) is a basis element in \so(2).
Fix \( n\in \N  \), \( n\geq 2 \), and set 
\begin{equation}
\label{classpsi}
\psi _{m}(z)=(z-\imath )^{m}(z+\imath )^{-m-n},\quad m=0,1,2,\dots .
\end{equation}
Then 
\begin{equation}
\varrho (E-F)\psi _{m}=\imath (2m+n)\, \psi _{m}
\end{equation}
and 
\begin{eqnarray}
 &  & \varrho (E)\psi _{m}=\frac{1}{2}\imath m\, 
\psi _{m-1}+\frac{1}{2}\imath (2m+n)\, \psi _{m}+\frac{1}{2}\imath (m+n)\, 
\psi _{m+1},\\
 &  & \varrho (F)\psi _{m}=\frac{1}{2}\imath m\, 
\psi _{m-1}-\frac{1}{2}\imath (2m+n)\, \psi _{m}+\frac{1}{2}\imath (m+n)\, 
\psi _{m+1},\\
 &  & \varrho (H)\psi _{m}=m\, \psi _{m-1}-(m+n)\, \psi _{m+1}.\label{classH} 
\end{eqnarray}
Note, however, that the eigen-functions \( \psi _{m} \) are not normalized.
The norm is 
\begin{equation}
\Vert \psi \Vert =2^{-n+1}\left( \pi \frac{m!\, (n-2)!}{(m+n-1)!}\right) ^{1/2}.
\end{equation}

\section{\protect\( q\protect \)-deformation of the discrete series}

Here we come to the main goal of the present paper, namely to a \( q \)-deformation
of the discrete series of representations of \( SL(2,\R ) \). Fix \( n\in \Z _{+} \)
and set \( \sigma =-n \) in (\ref{maction}). The relations (\ref{maction})
imply that
\begin{eqnarray}
  &  & K\cdot f(z)=q^{n/2}f(qz),\nonumber \\
  &  & E\cdot f(z)=\frac{z}{q-q^{-1}}(q^{n}K\cdot f(z)-K^{-1}\cdot f(z)),
  \label{qdifference} \\
  &  & F\cdot f(z)=
  -\frac{1}{(q-q^{-1})z}(q^{-n}K\cdot f(z)-K^{-1}\cdot f(z)).\nonumber 
\end{eqnarray}
Introduce the functions
\begin{equation}
  \label{qpsi}
  \psi _{m}(z)=\prod ^{m-1}_{j=0}(q^{2(j+n)}z-\imath )/
  \prod ^{m+n-1}_{j=0}(q^{2(j-m)}z+\imath ).
\end{equation}
Then it holds \-
\begin{equation}
  (q^{2n}EK-FK)\cdot \psi _{m}=\imath q^{n}[2m+n]_{q}\, \psi _{m}
\end{equation}
and 
\begin{eqnarray}
  EK^{-1}\cdot \psi _{m} & = & \imath 
  \frac{q^{4m+n-1}(1+q^{2})}{(1+q^{4m+2n-2})(1+q^{4m+2n})}[m]_{q^{2}}\, 
  \psi _{m-1}\nonumber \\
  &  & +\imath 
  \frac{q^{4m+2}(1+q^{2(n-1)})}{(1+q^{4m+2n-2})(1+q^{4m+2n+2})}[2m+n]_{q}\, 
  \psi _{m}\\
  &  & +\imath 
  \frac{q^{4m+n+1}(1+q^{2})}{(1+q^{4m+2n})(1+q^{4m+2n+2})}[m+n]_{q^{2}}\, 
  \psi _{m+1},\nonumber \\
  FK^{-1}\cdot \psi _{m} & = & \imath 
  \frac{q^{8m+5n-5}(1+q^{2})}{(1+q^{4m+2n-2})(1+q^{4m+2n})}[m]_{q^{2}}\, 
  \psi _{m-1}\nonumber \\
  &  & -\imath 
  \frac{q^{4m+2n}(1+q^{2(n-1)})}{(1+q^{4m+2n-2})(1+q^{4m+2n+2})}[2m+n]_{q}\, 
  \psi _{m}\\
  &  & +\imath 
  \frac{q^{n-3}(1+q^{2})}{(1+q^{4m+2n})(1+q^{4m+2n+2})}[m+n]_{q^{2}}\, 
  \psi _{m+1},\nonumber \\
  K^{-2}\cdot \psi _{m} & = & -(q^{2}-q^{-2})
  \frac{q^{6m+3n-2}}{(1+q^{4m+2n-2})(1+q^{4m+2n})}[m]_{q^{2}}\, 
  \psi _{m-1}\nonumber \\
  &  & +
  \frac{q^{4m+2}(1+q^{2(n-1)})(1+q^{2})}{(1+q^{4m+2n-2})(1+q^{4m+2n+2})}\, 
  \psi _{m}\label{qK} \\
  &  & +(q^{2}-q^{-2})
  \frac{q^{2m+n}}{(1+q^{4m+2n})(1+q^{4m+2n+2})}[m+n]_{q^{2}}\, 
  \psi _{m+1}.\nonumber 
\end{eqnarray}
Actually, relations (\ref{qpsi}-\ref{qK}) can be considered as 
\( q \)-deformations 
of relations (\ref{classpsi}-\ref{classH}). 



\begin{thebibliography}{**}
\bibitem{DM}Dobrev V.K., Moylan P.: 
Phys. Lett. B \textbf{315} (1993) 292-298. 
\bibitem{Woronowicz}Woronowicz S.L.: 
Lett. Math. Phys. \textbf{23} (1991) 251-263. 
\bibitem{Tarlini}Bonechi F., Ciccoli N., Giachetti R., Sorace E., Tarlini M.: 
Lett. Math. Phys. \textbf{49} (1999) 17-31. 
\bibitem{PS}\Stovicek P.: 
Lett. Math. Phys. \textbf{47} (1999) 125-138.
\bibitem{Joseph}Joseph A.: 
\emph{Quantum Groups and Their primitive Ideals}, Springer-Verlag,
Berlin, 1995. 
\end{thebibliography}
\end{document}